\RequirePackage{fix-cm}
\documentclass[smallextended]{svjour3}       
\smartqed  
\usepackage{graphicx}
\usepackage{amsfonts}
\usepackage{bm}
\usepackage{amsmath}
\begin{document}

\title{Improved Canonical Dual Algorithms for the Maxcut Problem
}


\author{Xiaojun Zhou         \and
        David Yang Gao       \and
        Chunhua Yang 
}


\institute{Xiaojun Zhou \at
              School of Science, Information Technology and Engineering, University of Ballarat, Victoria 3353, Australia.\\
              School of Information Science and Engineering, Central South University, Changsha 410083, China.\\
           \and
           David Yang Gao \at
              School of Science, Information Technology and Engineering, University of Ballarat, Victoria 3353, Australia.\\
           \and
           Chunhua Yang \at
              School of Information Science and Engineering, Central South University, Changsha 410083, China.
}

\date{Received: date / Accepted: date}

\maketitle

\begin{abstract}
By introducing a quadratic perturbation to the canonical dual of the maxcut problem, we transform the integer programming problem into a concave maximization problem over a convex positive domain under some circumstances, which can be solved easily by the well-developed optimization methods. Considering that there may exist no critical points in the dual feasible domain, a reduction technique is used gradually to guarantee the feasibility of the reduced solution, and a compensation technique is utilized to strengthen the robustness of the solution. The similar strategy is also applied to the maxcut problem with linear perturbation and its hybrid with quadratic perturbation. Experimental results demonstrate the effectiveness of the proposed algorithms when compared with other approaches.
\end{abstract}

\section{Introduction}
As one of the 21 NP-hard combinatorial optimization problems \cite{karp}, the maxcut problem has drawn a great deal of attention for several decades, due to its practical applications in circuit layout design, statistical physics, classification, and network analysis (for much more, the reader can refer to \cite{boros,poljak}).\\
\indent There are special classes of graphs, for instance, the planar graphs and graphs with no $K_5$ minor, for which the  maxcut problem can be solvable; however, it is impossible to find an algorithm in polynomial time due to its NP-hardness in general. Methods to solve the maxcut problem can be classified into two categories: direct and indirect. The direct method starts from a feasible solution, and then it ameliorates the solution iteratively with various strategies. Randomized heuristics which combine a greedy randomized adaptive search procedure, a variable neighborhood search, and a path-relinking intensification heuristic in \cite{festa}, advanced scatter search in \cite{marti}, a variant of spectral partitioning in \cite{trevisan}, and a new tabu search algorithm in \cite{kochenberger} belong to this kind. On the other hand, the indirect method focuses on the relaxation of the primal or the Lagrangian dual of the maxcut models, and then it uses some procedures to make the relaxation solution feasible to orginal problem, in which, relaxation means that relaxing some of the constraints and extending the objective function to larger space. Based on the semidefinite programming (SDP) relaxation, Goemans and Williamson \cite{goemans} proposed the well-known randomized approximation algorithm which can achieve a performance guarantee of  0.878. Based on the rank-two relaxation, a continuous optimization heuristic was constructed for approximating the maxcut problem \cite{burer}. Based on the Lagrangian dual relaxation, smoothing heuristics were presented in \cite{alperin}. In \cite{krishnan}, a polyhedral cut and price approach was investigated to solve the maxcut problem, in which, at the pricing phase, an interior point cutting plane algorithm was used to solve the Lagrangian dual of the SDP relaxation; at the cutting phase, cutting planes based on the polyhedral theory were added to the primal problem in order to improve the the SDP relaxation.  In \cite{xu}, a feasible direction method was designed to solve the continuous nonlinear programming problem by the relaxed maxcut model. In a branch and bound setting, a dynamic bundle method was used to solve the Lagrangian dual of the semidefinite maxcut relaxation\cite{rendl}.\\
\indent The canonical dual is a generalized Lagrangian dual, which originated in the late 1980s by Gao and Strang, and has developed significantly in recent years, which provides a new and potentially useful methodology for solving integer programming problems \cite{gao1,gao2}. As shown in \cite{fang,gao3}, through special canonical dual transformation, the integer programming problems can be converted to a concave maximization problem over a convex positive domain without duality gap, which can be solved easily by well-developed optimization algorithms under some circumstances. In \cite{wang}, a canonical dual approach was proposed for the maxcut problem by adding a linear perturbation term to the primal objective function; however, the strategy is not so robust to guarantee the solution to be located in the primal feasible domain. In this paper, we introduce a quadratic perturbation to the canonical dual of the maxcut problem. Furthermore, for some cases when there exist no critical points in the dual feasible space, a reduction technique is used gradually to guarantee the feasibility of the reduced solution and a compensation technique is utilized to strengthen the robustness of the solution. We also apply the strategy to the maxcut problem with linear perturbation and its hybrid with quadratic perturbation. Finally, experimental results are given to testify the effectiveness of the proposed algorithms.
\section{The Maxcut Problem}
Let $G = (V,E)$ be an undirected graph with edge weight $w_{ij}$ on $n+1 = |V|$ vertices and $m = |E|$ edges, for each edge $(i,j) \in E$, the maxcut problem is to find a
subset $S$ of the vertex set $V$ such that the total weight of the edges between $S$ and its complementary subset $\bar{S} = V \backslash S$ is as large as possible.
\subsection{LP based maxcut model}
Considering a variable $y_{ij}$ for each edge $(i,j) \in E$, and assuming $y_{ij}$ to be 1 if $(i,j)$ is in the cut, and 0 otherwise, the maxcut problem can be modeled as the following linear programming (LP) optimization problem:
\begin{eqnarray}\label{eq0}
    \max ~W(\mathbf{y}) & = & \sum_{i=1}^{n+1}\sum_{i<j,(i,j) \in E} w_{ij} y_{ij} \nonumber \\
         \mathrm{s.t.} ~ & & \mathbf{y} \mathrm{~is~the~incidence~vector~of~a~cut}.
\end{eqnarray}
Here the incidence vector $\mathbf{y} = \{y_{ij}\} \in \mathbb{R}^m$, where the $m$ is the number of edges in the graph.\\
\indent Let CUT($G$) denote the convex hull of the incidence vectors of cuts in $G$. Since maximizing a linear function over a set of points equals to maximizing it over the convex hull of this set of points, we can rewrite (\ref{eq0}) to the following form:
\begin{eqnarray}\label{eq01}
    \max ~W(\mathbf{y}) & = & \mathbf{c}^T \mathbf{y} \nonumber \\
         \mathrm{s.t.} & & \mathbf{y} \in \mathrm{CUT}(G).
\end{eqnarray}
where $c = \{w_{ij}\} \in \mathbb{R}^m$.
\subsection{SDP based maxcut model}
For a bipartition $(S,\bar{S})$, with $y_i = 1$ if $i \in S$, and $y_i = -1$ otherwise, the maxcut problem can be formulated as the following semidefinite programming optimization problem:
\begin{eqnarray}\label{eq1}
    \max ~W(\mathbf{y}) & = &\frac{1}{4} \sum_{i=1}^{n+1}\sum_{j=1}^{n+1} w_{ij}(1 - y_i y_j) \nonumber \\
    \mathrm{s.t.} ~ & & \mathbf{y} \in\{-1,1\}^{n+1}.
\end{eqnarray}
\indent Without loss of generality, if we fix the value of the last variable at 1, then the problem (\ref{eq1}) is equivalent to the following form (primal problem):
\begin{eqnarray}\label{eq2}
(\mathcal{P}) : \min\Big\{P(\mathbf{x}) = \frac{1}{2} \langle \mathbf{x}, Q \mathbf{x} \rangle - \langle \mathbf{x},\mathbf{c} \rangle: \mathbf{x} \in \{-1,1\}^{n}\Big\},
\end{eqnarray}
where, $\langle \mathbf{u}, \bar{\mathbf{u}} \rangle = \mathbf{u}^T \bar{\mathbf{u}}$ denotes the bilinear form,  $Q$ is a symmetric matrix with $Q_{ij} = w_{ij}(i,j = 1,2,\cdots, n)$, and $\mathbf{c} = -(w_{1(n+1)}, \cdots, w_{n(n+1)})^T$.  It is not difficult to find that a optimal solution $\mathbf{x}^{*}$ to problem (\ref{eq2}) corresponds to a optimal solution $(\mathbf{x}^{*},1)$ of original problem (\ref{eq1}).\\
\indent The canonical dual algorithms in this paper are based on the SDP model.
\section{Primal problem with linear perturbation}
In \cite{wang}, a linear non-zero vector $\triangle \mathbf{c} = (\triangle c_1, \cdots, \triangle c_n) \in \mathbb{R}^n$ is added to the primal problem and then we can get the linearly perturbed problem
\begin{eqnarray}
(\mathcal{P}_l) : \min\Big\{P(\mathbf{x}) = \frac{1}{2} \langle \mathbf{x}, Q \mathbf{x} \rangle - \langle \mathbf{x},\mathbf{c} + \triangle \mathbf{c} \rangle: \mathbf{x} \in \{-1,1\}^{n}\Big\}.
\end{eqnarray}
\indent It proved that, if $\sum_{i=1}^n|\triangle c_i| < 1$, the perturbed optimal solution is equivalent to the solution of primal problem, and if $\sum_{i=1}^n|\triangle c_i| = \alpha W(\widetilde{\mathbf{x}})$, where $\widetilde{\mathbf{x}}$  is the optimal
solution of the perturbed problem and $\alpha > 0$ is a constant, then $W(\widetilde{\mathbf{x}}) \geq \frac{1}{1 + \alpha} W^{*}$, here, $W^{*}$ is the optimal value of the primal problem.\\
\indent The problem of this linear perturbation is that, for one thing, it is not easy to find such an appropriate perturbation, and for another, if $\sum_{i=1}^n|\triangle c_i| > 1$, only approximate solution can be obtained.
\section{Canonical dual problem with quadratic perturbation}
We give the notations $\mathbf{s} \circ \mathbf{t}:= [s_1t_1, \cdots, s_mt_m]$ and $\mathbf{s} \oslash \mathbf{t}:= [ s_1/t_1, \cdots, s_m/t_m]$ to denote the Hadamard product and quotient for any two vector $\mathbf{s}, \mathbf{t} \in \mathbb{R}^m$, $Diag(\bm\nu)$ to represent a diagonal matrix with components of the vector $\bm\nu$ as its elements, and $\mathbf{e}$ to stand for a vector of all ones.\\
\indent By the fact $\mathbf{x} \circ \mathbf{x} = \mathbf{e}$ and $[\mathbf{x} \circ \mathbf{x}]^T \bm\alpha = \mathbf{x}^T Diag(\bm\alpha)\mathbf{x} = \mathbf{e}^T \bm\alpha$, then the problem $(\mathcal{P})$ is identical to the following $\alpha-$perturbed problem:
\begin{equation}
(\mathcal{P}_{\alpha}): \min \Big\{  P_{\alpha}(\mathbf{x}) = \frac{1}{2} \langle \mathbf{x},(Q+Diag(\bm\alpha)) \mathbf{x} \rangle - \langle \mathbf{x},\mathbf{c} \rangle - d_{\alpha}|\mathbf{x} \in \{-1,1\}^n \Big\},
\end{equation}
where $d_{\alpha} = \frac{1}{2}\langle \mathbf{e},\bm\alpha \rangle$, $\bm\alpha \in \mathbb{R}^n$ is a parametrical vector.\\
\indent Introducing a quadratic geometrical operator
\begin{eqnarray}
\bm\epsilon = \mathrm{\Lambda}(\mathbf{x}) = \frac{1}{2} (\mathbf{x} \circ \mathbf{x} - \mathbf{e})
\end{eqnarray}
and following the standard procedures of canonical dual methodology, the $\alpha-$perturbed canonical dual problem is obtained as follows
\begin{eqnarray}
(\mathcal{P}_{\alpha}^d): \max \Big\{P_{\alpha}^d(\bm{\sigma}) = -\frac{1}{2} \langle G_{\alpha}^{-1}(\bm\sigma)\mathbf{c},\mathbf{c}\rangle - \frac{1}{2} \langle \mathbf{e},\bm\sigma \rangle - d_{\alpha}|\bm\sigma \in \mathcal{S}^{+}_{\alpha}\Big\},
\end{eqnarray}
where, $G_{\alpha}(\bm\sigma) = Q + Diag(\bm\alpha + \bm\sigma)$, and
\begin{eqnarray}
\mathcal{S}^{+}_{\alpha}= \{\bm\sigma \in \mathbb{R}^n |G_{\alpha}(\bm\sigma) \succ 0 \}.
\end{eqnarray}
\indent Further more, a quadratic penalty term is added to the $\alpha-$perturbed canonical dual problem, and then we get the $\beta-$ perturbed canonical dual problem
\begin{equation}\label{eq7}
(\mathcal{P}_{\alpha \beta}^d): \max \Big\{P_{\alpha \beta}^d(\bm{\sigma}) = -\frac{1}{2} \langle G_{\alpha}^{-1}(\bm\sigma)\mathbf{c},\mathbf{c}\rangle - \frac{1}{2}\sum\limits_{i=1}^{n}(\frac{\sigma_i^2}{\beta_i} + \sigma_i) - d_{\alpha}|\bm\sigma \in \mathcal{S}^{+}_{\alpha}\Big\}.
\end{equation}
\begin{theorem}
Suppose that for a given perturbation vector $\bm\alpha$ such that $Q + Diag(\bm\alpha) \prec 0$, and $\beta_i \gg 0(i=1,2,\cdots,n)$ is big enough, if $\bar{\bm\sigma} \in \mathcal{S}^{+}_{\alpha}$ is a critical point of the $\beta-$ perturbed canonical dual problem, then
\begin{eqnarray}\label{eq8}
\bar{\mathbf{x}} = G_{\alpha}^{-1}(\bar{\bm\sigma})\mathbf{c}, \mathbf{x}^{*} = round(\bar{\mathbf{x}}),
\end{eqnarray}
$\mathbf{x}^{*}$ is a global solution to the $(\mathcal{P})$.
\end{theorem}
\proof
First, we prove that $\mathbf{x}^{*}$ defined by (\ref{eq8}) is feasible of $(\mathcal{P})$.\\
\indent Considering that
\begin{eqnarray*}
\frac{\partial P_{\alpha \beta}^d(\bm{\sigma})}{\partial \bar\sigma_i} & = & -\frac{1}{2}\frac{\partial [\mathbf{c}^T G_{\alpha}^{-1}(\bm\sigma) \mathbf{c}]}{\partial \bar\sigma_i} - \frac{\bar\sigma_i}{\beta_i} - \frac{1}{2} \\
& = & \frac{1}{2} \mathbf{c}^T G_{\alpha}^{-1}(\bar{\bm\sigma}) \frac{ \partial G_{\alpha}(\bm\sigma)} {\partial \bar\sigma_i} G_{\alpha}^{-1}(\bar{\bm\sigma}) \mathbf{c} - \frac{\bar\sigma_i}{\beta_i} - \frac{1}{2} \\
& = & \frac{1}{2} (\bar{\mathbf{x}} \circ \bar{\mathbf{x}})_i - \frac{\bar\sigma_i}{\beta_i} - \frac{1}{2} \\
& = &  \frac{1}{2}  \bar{x}^2_i - \frac{\bar\sigma_i}{\beta_i} - \frac{1}{2},
\end{eqnarray*}
if $\bar{\bm\sigma}$ is a critical point of $(\mathcal{P}_{\alpha \beta}^d)$, then $\frac{1}{2}  \bar{x}^2_i - \frac{\bar\sigma_i}{\beta_i} - \frac{1}{2} = 0$, namely, $\bar{x}^2_i = \frac{2\bar\sigma_i}{\beta_i} + 1$. \\
\indent If $\beta_i \gg 0$ is big enough, it is easy to find that $x^{*}_i = round(\bar{x}_i) = \pm 1$, which is a feasible solution to $(\mathcal{P})$. \\
\indent Then, we prove that $\mathbf{x}^{*}$ is a global solution.\\
\indent By the fact that $Q + Diag(\bm\alpha) \prec 0$, we can conclude that $\bar{\bm\sigma} > 0$ due to $G_{\alpha}(\bar{\bm\sigma}) \succ 0$. Therefore, the $\beta-$ perturbed canonical dual problem $(\mathcal{P}_{\alpha \beta}^d)$ is strictly concave on $\mathcal{S}^{+}_{\alpha}$, and if the critical point $\bar{\bm\sigma} \in \mathcal{S}^{+}_{\alpha}$, it must be a unique maximizer of problem $(\mathcal{P}_{\alpha \beta}^d)$. Consequently, the corresponding $\mathbf{x}^{*}$ defined by (\ref{eq8}) must be unique. On the other hand, we have
\begin{eqnarray*}
P_{\alpha}(\bar{\mathbf{x}}) &=& \frac{1}{2} \langle \bar{\mathbf{x}},(Q+Diag(\bm\alpha)) \bar{\mathbf{x}} \rangle - \langle \bar{\mathbf{x}},\mathbf{c} \rangle - d_{\alpha} \\
&=& \frac{1}{2} \langle \bar{\mathbf{x}}, G_{\alpha}(\bar{\bm\sigma}) \bar{\mathbf{x}} \rangle - \langle \bar{\mathbf{x}},\mathbf{c} \rangle - \frac{1}{2}\langle \bar{\bm\sigma},  \mathbf{e} \rangle - d_{\alpha}\\
&=&\min_{\mathbf{x} \in \mathbb{R}^n} \frac{1}{2} \langle \mathbf{x}, G_{\alpha}(\bar{\bm\sigma}) \mathbf{x} \rangle - \langle \mathbf{x},\mathbf{c} \rangle - \frac{1}{2}\langle \bar{\bm\sigma},  \mathbf{e} \rangle - d_{\alpha}\\
&\leq& \min_{\mathbf{x} \in \{-1,1\}^n} \frac{1}{2} \langle \mathbf{x},(Q+Diag(\bm\alpha)) \mathbf{x} \rangle - \langle \mathbf{x},\mathbf{c} \rangle - d_{\alpha},
\end{eqnarray*}
that is to say, $\mathbf{x}^{*} = round(\bar{\mathbf{x}})$ is the global solution of $(\mathcal{P}_{\alpha})$ and so the global solution of $(\mathcal{P})$.
\section{Canonical dual algorithms}
Now, we present a gradient-based iterative method for the $\beta-$ perturbed canonical dual problem.\\
\indent It is not difficult to find that the search direction (negative gradient) of the negative canonical dual function $-P_{\alpha \beta}^d(\bm{\sigma})$ is
\begin{eqnarray}
d = \frac{1}{2} \mathbf{x} \circ \mathbf{x} - \bm{\sigma} \oslash \bm \beta - \frac{1}{2} \mathbf{e},
\end{eqnarray}
and step size of each component of search direction is obtained by a golden section search technique from $[0,\alpha_{\mathrm{max}}]$, where $\alpha_{\mathrm{max}}$ is determined by making the corresponding diagonal component of
$G_{\alpha}(\bm\sigma)$ positive. \\
\indent We use three types of termination criteria to stop the procedure. The details of the algorithm is give in the following
\begin{table}[!htbp]
\begin{tabular}{l}
  \hline
  Algorithm 1\\
  \hline
  1: Initialization, set $\bm\sigma_k = \bm\sigma^0$, $\bm\alpha$, $\bm\beta$, $\epsilon$ and $k = 0$\\
  2: \textbf{while} (1) \textbf{do}\\
  3: ~~$\mathbf{x}_k \leftarrow (Q+Diag(\bm\alpha + \bm\sigma_k))^{-1}\mathbf{c}$, $\mathbf{d}_k \leftarrow \frac{1}{2} \mathbf{x}_k \circ \mathbf{x}_k - \bm{\sigma}_k \oslash \bm \beta - \frac{1}{2} \mathbf{e}$\\
  4: ~~$\bm\sigma_{k+1} = \bm\sigma_k + a_k \times \mathbf{d}_k$  \\
  5: ~~$\mathbf{x}_{k+1} \leftarrow (Q+Diag(\bm\alpha + \bm\sigma_{k+1}))^{-1}\mathbf{c}$\\
  6: ~~\textbf{if} $\|d_k\| \leq \varepsilon ~ \mathrm{or} ~ \|\frac{\mathbf{x}_{k+1} - \mathbf{x}_{k}}{\mathbf{x}_{k}}\| \leq \varepsilon ~ \mathrm{or} ~ \| \frac{\bm\sigma_{k+1} - \bm\sigma_{k}}{\bm\sigma_{k}} \| \leq \varepsilon$ \textbf{then}\\
  7: ~~~~ $\mathbf{x}^{*} \leftarrow round(\mathbf{x}_{k+1})$, \textbf{break} \\
  8: ~~\textbf{end if}\\
  9: ~~$k \leftarrow k + 1$ \\
  10: \textbf{end while}\\
  \hline
\end{tabular}
\end{table}
\\
where, $\varepsilon$ is the numerical precision, and the initial point $\bm\sigma^0$ should make $G_{\alpha}(\bm\sigma^0)$ positive.
\\
\textbf{Remark 1} If the final critical point $\bm\sigma^{*}$ of $(\mathcal{P}_{\alpha \beta}^d)$ is in $\mathcal{S}^{+}_{\alpha}$, according to Theorem 1, we can conclude that the corresponding $\mathbf{x}^{*}$ is the global solution to $(\mathcal{P})$.\\
\indent However, the critical point $\bm\sigma^{*}$ is not always in $\mathcal{S}^{+}_{\alpha}$, or in $\mathcal{S}^{+}_{\alpha}$, we can not find a critical point $\bm\sigma^{*}$. In this case, we can find that some components of $\mathbf{x}^{*}$ are feasible, and other components are infeasible. Thus, there must exist a transformation $N \in \mathbb{R}^{n \times m}$, such that
\begin{eqnarray}
\mathbf{x}^{*} = \mathbf{x}_p + N \mathbf{x}_h,
\end{eqnarray}
where, $\mathbf{x}_p \in \mathbb{R}^n$ is the particular solution with infeasible components zeros, and $\mathbf{x}_h \in \mathbb{R}^m$ is called the reduced solution.\\
\indent Considering that
\begin{eqnarray}
P(\mathbf{x}^{*}) & = & P(\mathbf{x}_p + N \mathbf{x}_h) \nonumber \\
& = & \frac{1}{2}(\mathbf{x}_p + N \mathbf{x}_h)^T Q (\mathbf{x}_p + N \mathbf{x}_h) - c^T(\mathbf{x}_p + N \mathbf{x}_h) \nonumber \\
& = & \frac{1}{2} \mathbf{x}_h^T (N^T Q N) \mathbf{x}_h - \mathbf{x}_h^T(N^T \mathbf{c} - N^T Q \mathbf{x}_p) + \frac{1}{2}\mathbf{x}_p^T Q \mathbf{x}_p - c^T\mathbf{x}_p,
\end{eqnarray}
as a result, we have to solve the reduced problem in further
\begin{eqnarray}
P(\mathbf{x}_h) & = & \frac{1}{2} \langle \mathbf{x}_h, Q_h \mathbf{x}_h \rangle - \langle \mathbf{x}_h,\mathbf{c}_h \rangle,
\end{eqnarray}
where, $Q_h = N^T Q N, c_h = N^T \mathbf{c} - N^T Q \mathbf{x}_p$.\\
\indent The reduction technique above can be regarded as ``greedy criterion", which belongs to local search to some extent. To compensate for the lost information in the reduction process,  we use a simple compensation technique, that is, replacing each binary component in $\mathbf{x}^{*}$ with its counterpart, to strength the robustness of the solution. The details of the improved canonical dual algorithm are given below:
\\
\begin{table}[!htbp]
\begin{tabular}{l}
  \hline
  Algorithm CDA1\\
  \hline
  1: Using Algorithm 1 as a subroutine to get $\mathbf{x}^{*}$\\
  2: \textbf{if} the corresponding $\bm\sigma^{*} \not\in \mathcal{S}^{+}_{\alpha}$ \textbf{then} \\
  3: ~~obtaining current $\mathbf{x}_h$ by the reduction technique \\
  4: ~~\textbf{while} $length(\mathbf{x}_h) \geq 1$ \textbf{do} \\
  5: ~~~~~solving the current reduced problem $P(\mathbf{x}_h)$ by Algorithm 1 again\\
  6: ~~~~~\textbf{if} the corresponding $\bm\sigma^{*}_h \in \mathcal{S}^{+}_{\alpha(h)}$ \textbf{then}\\
  7: ~~~~~~~$\mathbf{x}^{*} = \mathbf{x}_p + N \mathbf{x}_h$, \textbf{break}\\
  8: ~~~~~\textbf{else}\\
  9: ~~~~~~~obtaining current $\mathbf{x}_h$ by the reduction technique again\\
  10: ~~~~\textbf{end if}\\
  11: ~~\textbf{end while}\\
  12: \textbf{end if}\\
  13: Using the compensation technique to strength the robustness of the solution\\
  \hline
\end{tabular}
\end{table}
\\
We also apply the similar reduction technique and compensation technique to the maxcut problem with linear perturbation in \cite{wang} and its hybrid with quadratic perturbation in this paper, so we get the corresponding algorithms named CDA2 and CDA3.
\section{Experimental results}
To testify the effectiveness of the proposed algorithms, firstly, we give a random example to show that in some case, the Algorithm 1 can achieve a global solution without additional strategies. Then we compared our algorithms with other approaches for some published TSPLIB (see in \cite{reinelt}) instances results. Finally, some large instances in TSPLIB are tested to show the capacity of our algorithms. We run the proposed procedures in MATLAB R2010b on Intel(R) Core(TM) i3-2310M CPU @2.10GHz under Window 7 environment.
\\
\indent \textbf{Example 1} Consider the following 10-dimensional problem(n=9)  with randomly selected matrix $Q$ and $\mathbf{c}$, $\bm\alpha$ and $\bm\beta$
\begin{eqnarray*}
Q =
\begin{pmatrix}
     0  &   6   &   4   &  8   &   4   &   5    &  5    &  6    &  8\\
     6  &   0   &   3   &  9   &  4    &  5     & 5    &  8   &   7\\
     4  &   3   &   0   &  6   &   2   &   4  &    7   &   5  &    4\\
     8  &   9   &   6   &  0   &  7    &  4   &   7   &   6   &   6\\
     4  &   4   &   2   &  7   &   0    &  7   &   5    &  4    &  6\\
     5  &   5   &   4   &  4   &    7   &   0   &   0    &  2   &   7\\
     5  &   5   &   7   &  7   &   5   &   0   &   0  &    4   &   5\\
     6  &   8   &   5   &  6   &  4   &   2    &  4    &  0    &  2\\
     8  &   7   &   4   &  6   &  6   &   7  &    5   &   2    &  0\\
\end{pmatrix},
\mathbf{c} =
\begin{pmatrix}
 2\\
 5\\
 3\\
 5\\
 2\\
 5\\
 4\\
 5\\
 7\\
\end{pmatrix},
\bm\alpha =
\begin{pmatrix}
  -17.3208\\
   -2.8050\\
  -36.5410\\
   -1.1174\\
  -38.3706\\
  -77.0470\\
  -20.1651\\
  -31.3002\\
  -34.9461\\
\end{pmatrix},
\bm\beta =
\begin{pmatrix}
  605.7162\\
  601.1675\\
  330.2360\\
  277.4284\\
  674.9582\\
  540.0750\\
  537.7345\\
  690.3018\\
  371.8627\\
\end{pmatrix}
\end{eqnarray*}
we can get the dual critical solution \\
$\bm\sigma^{*} = (34.0286,19.9327,50.1747,12.6699,55.9428,88.7105,30.2908,45.0242,45.6742)^T$, \\
and the corresponding\\
$\mathbf{x}^{*} = round(G^{-1}(\bar{\bm\sigma})c) = (1,1,1,-1,1,-1,-1,-1,-1)$.\\
We can check that $\bm\sigma^{*}$ is in the dual feasible domain $\mathcal{S}^{+}_{\alpha}$, that is to say, the corresponding $\mathbf{x}^{*}$ is the global solution according to
Theorem 1. \\
\indent \textbf{Example 2} We compare our algorithms with the published results found in \cite{goemans}(\textit{GW}) and \cite{burer,commander}(\textit{circut}). From Table \ref{tab:1}, we find that all of the canonical dual algorithms can achieve the same or better solutions than that of (\textit{circut}) and (\textit{GW}), with comparable shorter time.
\begin{table}[!htbp]
\caption{Comparison with other algorithms for medium-size TSPLIB instances}
\label{tab:1}
\begin{tabular}{c|c|c|c|c|c|c|c|c|c|c}
  \hline
  Instances  &  \textit{circut} & time &  \textit{GW} & time  & \textit{CDA1} & time& \textit{CDA2} & time& \textit{CDA3} & time\\
  \hline
  burma14   & 283     & 0.046 & -       & -     & 283    & 0.216 & 283    & 0.084 &  283 & 0.095 \\
  gr17      & 24986   & 0.047 & -       & -     & 24986  & 0.087 & 24986  & 0.272 & 24986& 0.189\\
  bays29    & 53990   & 1.109 & -       & -     & 53990  & 0.349 & 53990  & 0.159 & 53990& 0.152\\
  dantzig42 & 42638   & 1.75  & 42638   & 43.35 & 42638  & 0.540 & 42638  & 0.204 & 42638& 0.203\\
  gr48      & 320277  & 3.672 & 320277  & 26.17 & 320277 & 0.346 & 320277 & 0.287 & 320277&0.294\\
  hk48      & 771712  & 2.516 & 771712  & 66.52 & 771712 & 1.368 & 771712 & 0.426 & 771712&0.510\\
  gr96      & 105328  & 14.250& 105295  & 531.50& 105328 & 0.550 & 105328 & 0.975 & 105328&0.693\\
  kroA100   & 5897368 & 2.359 & 5897392 & 420.83& 5897392& 1.444 & 5897392& 1.053 & 5897392&0.785\\
  kroB100   & 5763020 & 2.531 & 5763047 & 917.47& 5763047& 1.440 & 5763047& 1.021 & 5763047&0.890\\
  kroC100   & 5890745 & 2.500 & 5890760 & 398.78& 5890760& 1.258 & 5890760& 0.986 & 5890760&0.898\\
  kroD100   & 5463250 & 2.547 & 5463250 & 469.48& 5463250& 1.571 & 5463250& 1.115 & 5463250&0.843\\
  kroE100   & 5986587 & 2.500 & 5986591 & 375.68& 5986591& 1.476 & 5986591& 1.476 & 5986591&0.888\\
  gr120     &   -     &   -   & 2156667 & 754.87& 2156667& 2.752 & 2156667& 1.978 & 2156667&1.358\\
  \hline
\end{tabular}
\end{table}
\\
\indent \textbf{Example 3} Other instances from TSPLIB are used to test the capacity of the proposed algorithms, and we test the size of the instances up to 500 with no more than 2 minutes. In the meanwhile, we use a toolbox called YALMIP \cite{lofberg} to solve the same problems for comparison. The results also show the different performance of the three canonical dual algorithms (CDA3$ >$ CDA1$ >$ CDA2), which indicate that the hybrid with linear and quadratic perturbation is a better choice, especially for gil262.
\begin{table}[!htbp]
\caption{Comparison with YALMIP for large-size TSPLIB instances}
\label{tab:2}
\begin{tabular}{c|c|c|c|c|c|c|c|c}
  \hline
  Instances & \textit{YALMIP} & time  & \textit{CDA1}  & time& \textit{CDA2}  & time&\textit{CDA3}  & time\\
  \hline
  ch130 & 1888108 & 3.448 & 1888108    & 3.036 & 1888108  & 2.608& 1888108   & 1.461\\
  ch150 & 2525606 & 4.007 & 2525626   & 3.767 & 2525626  & 3.400& 2525626   & 1.671\\
  d198  & 12938532 &  6.910 &12938532  & 9.256 & 12938532 & 16.725& 12938532 & 3.790\\
  gr202 & 195433 &  6.980 & 197098    & 7.600 & 197098   & 8.196& 197098    & 4.866\\
  gr229 &  1203249 & 9.463 & 1205180   & 9.513 & 1205180  & 7.697& 1205180   & 8.576\\
  gil262&  2131227 & 13.360 & 2149623    & 10.479& 2147693  & 10.352&2152173$^{*}$   & 6.892\\
  pr299 &  80047271 & 20.156 & 80324674   & 34.453& 80324674 & 19.710&80324674  & 14.165\\
  lin318&  59504077 & 24.221 & 59547803  & 17.888& 59547803 & 33.739&59547803  & 13.158\\
  fl417 & 76742924 & 61.888 & 76776296   & 28.738& 76776296 & 40.181&76776296  & 22.590\\
  pr439 & 277700818 &  73.775 & 278552199 & 46.210& 278552199& 55.039&278552199 & 37.913 \\
  d493  &  75651106 & 99.701 & 75740087  & 53.323& 75740087 & 66.291&75740087  & 45.986\\
  att532&  281648797 & 127.485 & 287417240 & 57.681& 287417240& 109.979&287417240 & 52.192\\
  \hline
\end{tabular}
\end{table}
\section{Conclusion}
In this paper, we have studied the maxcut problem with a quadratic perturbation term added to the canonical dual problem.
A gradient-based algorithm is designed to solve the quadratic perturbed canonical dual problem.
For some cases with no critical points in the dual feasible domain, a reduction technique and a compensation technique are introduced to improve the quality of the solution. The similar techniques are also applied to maxcut model with linear perturbation and its hybrid with
quadratic perturbation. Numerical results show that the proposed algorithms are effective in terms of both time complexity and the solution quality, in the same time, it indicates that the hybrid perturbation is a more promising approach.


\begin{thebibliography}{}
%
%
\bibitem{karp} Karp, R.M.: Reducibility among combinatorial problems, In: Complexity of Computer Computations, edited by R.E. Miller and J.W. Thatcher, Plenum Press, New York.  pp. 85--103 (1972)
\bibitem{boros} Boros, E., Hammer, P.: The max-cut problem and quadratic 0-1 optimization; polyhedral aspects, relaxations and bounds, Annals of Operations Research. 33, 151--180 (1991)
\bibitem{poljak} Poljak, S., Tuza, Z.: Maximum cuts and large bipartite subgraphs. In W. Cook, L. Lov\'{a}sz, and P. Seymour, editors, Combinatorial Optimization, DIMACS Series in Discrete Mathematics and Theoretical Computer Science, American Mathematical Society. 20, 181--244 (1995)
\bibitem{festa} Festa, P., Pardalos, P.M., Resende, M.G.C., Ribeiro, C.C.: randomized heuristics for the max-cut problem. Optimization Methods and Software. 7, 1033--1058 (2002)
\bibitem{marti} Mart\'{\i}, R., Duarte, A., Laguna, M.: Advanced scatter search for the max-cut problem. NFORMS Journal on Computing. 21(1), 26--38 (2009)
\bibitem{trevisan} Trevisan, L.: Max cut and the smallest eigenvalue. in STOC. 263--272 (2009)
\bibitem{kochenberger} Kochenberger, G.A., Hao, J.K., L\"{u}, Z.P., Wang, H.B., Glover, F.: Solving large scale max cut problems via tabu search. online in Journal of Heuristics. (2012)
\bibitem{goemans} Goemans, M.X., Williamson, D.P.: Improved approximation algorithms for maximum cut and satisfiability problems using semidifinite programming. Journal of the Association for Computing Machinery. 42(6), 1115--1145 (1995)
\bibitem{burer} Burer, S., Monteiro, R.D.C., Zhang, Y.: Rank-two relaxation heuristics for max-cut and other binary quadratic programs. SIAM Journal on Optimization. 12, 503--521 (2001)
\bibitem{alperin} Alperin, H., Nowak, I.: Lagrangian smoothing heuristics for max-cut. Journal of Heuristics. 11, 447--463 (2005)
\bibitem{krishnan}Krishnan,  K., Mitchell, J.E.: A semidefinite programming based polyhedral cut and price approach for the maxcut problem. Computational Optimization and Applications. 35, 51--71 (2006)
\bibitem{xu} Xu, C.X., He X.L. and Xu, F.M.: An effective continuous algorithm for approximate solutions of large scale max-cut problem. Journal of Computational Mathematics. 24(6), 749--760 (2006)
\bibitem{rendl} Rendl, F., Rinaldi, G., Wiegele,  A.: Solving max-cut to optimality by intersecting semidefinite and polyhedral relaxations.  Mathematical Programming. 121(2), 307--335 (2010)
\bibitem{gao1} Gao, D.Y., Strang, G.: Geometric nonlinearity: potential energy, complementary energy, and the gap function. Quarterly Journal of Applied Mathematics. (XLVII)(3), 487--504 (1989)
\bibitem{gao2} Gao, D.Y.: Caonical duality theory: Unified understanding and generalized solution for global optimization problems. Computers and Chemical Engineering. 33, 1964--1972 (2009)
\bibitem{fang} Fang, S.C., Gao, D.Y., Sheu, R.L., Wu, S.Y.: Canonical dual approach to solving 0-1 quadratic programming problems. Journal of Industrial and Management Optimization. 4(1), 125--142 (2008)
\bibitem{gao3} Gao, D.Y., Ruan, N.: Solutions to quadratic minimization problems with box and integer constraints. Journal of Global Optimization. 47, 463--484 (2010)
\bibitem{wang} Wang, Z.B., Fang, S.C., Gao, D.Y., Xing, W.X.: Canonical dual approach to solving the maximum cut problem. oneline in Journal of Global Optimization. (2012)
\bibitem{commander} Commander, C.W.: Maximum cut problem, max-cut. In: C. A. Floudas and P. M. Pardalos, editors, Encyclopedia of Optimization, Springer US. pp. 1991--1999 (2009)
\bibitem{reinelt} Reinelt,  G.: TSPLIB--A traveling salesman problem library. ORSA Journal on Computing. 3(4), 376--384 (1991)
\bibitem{lofberg} L\"{o}fberg, J.: YALMIP : A Toolbox for Modeling and Optimization in MATLAB. in Proceedings of the CACSD Conference, Taipei, Taiwan. (2004)

\end{thebibliography}


\end{document}